\documentclass{article}

\usepackage[active]{srcltx}
\usepackage{graphicx}
\usepackage{latexsym}
\usepackage[utf8]{inputenc}
\usepackage{amsfonts, amssymb}
\usepackage{amsmath}
\usepackage{enumerate}
\usepackage{graphicx}
\usepackage[labelsep=period,font=small]{caption}
\usepackage{mathrsfs}
\usepackage{authblk}
\usepackage{color}

\newtheorem{thm}{Theorem}[section]
\newtheorem{defn}[thm]{Definition}
\newtheorem{lemma}[thm]{Lemma}
\newtheorem{prop}[thm]{Proposition}

\newtheorem{rem}[thm]{Remark}
\newtheorem{cor}[thm]{Corollary}
\newtheorem{example}[thm]{Example}

\newtheorem{defthm}[thm]{Definition and Theorem}
\newtheorem{question}[thm]{Question}

\newcommand{\T}{{\mathcal T}}
\newcommand{\C}{{\mathbb C}}
\newcommand{\Chat}{{\widehat{\mathbb C}}}

\newcommand{\Om}{\Omega}
\newcommand{\om}{\omega}
\newcommand{\sm}{\setminus}
\newcommand{\ga}{\gamma}

\newcommand{\si}{\sigma}

\newcommand{\Laplace}{\Delta}

\newcommand{\eps}{\epsilon}
\newcommand{\N}{\mathbb N}
\newcommand{\D}{\mathbb D}

\newcommand{\KK}{\mathcal K}
\newcommand{\R}{{\mathbb R}}

\newcommand{\Sm}{\setminus}
\newcommand{\wtnu}{{\widetilde\nu}}

\newcommand{\dH}{\ensuremath{{\operatorname{d_H}}}}
\renewcommand{\DH}{\ensuremath{{\operatorname{D_H}}}}
\renewcommand{\d}{\ensuremath{\operatorname{d}}}

\newcommand{\diam}{\ensuremath{\operatorname{diam}}}
\newcommand{\e}{\ensuremath{\operatorname{e}}}
\newcommand{\dnu}{\ensuremath{\operatorname{d\nu}}}

\newcommand{\mapfromto}[3]{\hbox{\ensuremath{#1 : #2 \longrightarrow #3}}}

\DeclareMathOperator{\Cpct}{Cap}
\DeclareMathOperator{\Co}{Co}
\DeclareMathOperator{\Po}{Po}

\newcommand{\ALIGN}{\begin{align*}}
\newcommand{\ENDALIGN}{\end{align*}}
\newcommand{\ENUM}{\begin{enumerate}}
\newcommand{\ENUMa}{\begin{enumerate}[a.]}
\newcommand{\ENUMA}{\begin{enumerate}[A.]}
\newcommand{\ENUMi}{\begin{enumerate}[i)]}
\newcommand{\ENDENUM}{\end{enumerate}}
\newcommand{\ITMZ}{\begin{itemize}}
\newcommand{\ENDITMZ}{\end{itemize}}
\newcommand{\EQN}[1] { \begin{equation}\label{#1} }
\newcommand{\ENDEQN}{\end{equation}}
\newcommand{\THM}{\begin{thm}}
\newcommand{\EXA}{ \begin{example}}
\newcommand{\REFEXA}[1] { \begin{example}\label{#1} }
\newcommand{\ENDEXA}{\end{example}}
\newcommand{\REM}{ \begin{rem}}
\newcommand{\REFREM}[1] { \begin{rem}\label{#1}}
\newcommand{\ENDREM}{\end{rem}}
\newcommand{\REFTHM}[1] { \begin{thm}\label{#1} }
\newcommand{\RTHM}[1] { \begin{thm}[#1] }
\newcommand{\RREFTHM}[2] { \begin{thm}[#1]\label{#2} }
\newcommand{\ENDTHM}{\end{thm}}
\newcommand{\REFPROP}[1]{\begin{prop}\label{#1} }
\newcommand{\RREFPROP}[2]{\begin{prop}[#1]\label{#2} }
\newcommand{\RPROP}[1]{\begin{prop}[#1] }
\newcommand{\PROP}{\begin{prop}}
\newcommand{\ENDPROP}{\end{prop} }
\newcommand{\REFDEF}[1]{\begin{defn}\label{#1} }
\newcommand{\RREFDEF}[2]{\begin{defn}[#1]\label{#2} }
\newcommand{\DEF}{\begin{defn}}
\newcommand{\RDEF}[1] {\begin{defn}[#1]}
\newcommand{\ENDDEF}{\end{defn} }
\newcommand{\REFLEM}[1]{\begin{lemma}\label{#1} }
\newcommand{\RLEM}[1]{\begin{lemma}[#1] }
\newcommand{\RREFLEM}[2]{\begin{lemma}[#1]\label{#2} }
\newcommand{\LEM}{\begin{lemma}}
\newcommand{\ENDLEM}{\end{lemma} }
\newcommand{\CENTER}{\begin{center}}
\newcommand{\ENDCENTER}{\end{center}}
\newcommand{\REFCOR}[1]{\begin{cor}\label{#1} }
\newcommand{\RCOR}[1] {\begin{cor}[#1]}
\newcommand{\COR}{\begin{cor}}
\newcommand{\ENDCOR}{\end{cor}}
\newcommand{\RMRK}{\begin{remark}}
\newcommand{\ENDRMRK}{\end{remark}}
\newcommand{\REFDEFTHM}[1] { \begin{defthm}\label{#1} }
\newcommand{\RREFDEFTHM}[2] { \begin{defthm}[#1]\label{#2} }
\newcommand{\ENDDEFTHM}{\end{defthm}}
\newcommand{\QUE}{ \begin{question}}
\newcommand{\ENDQUE}{\end{question}}
\newcommand{\ENDCD}{\end{CD}}
\newcommand{\corref}[1]{Corollary~\ref{#1}}

\newcommand{\exaref}[1]{Example~\ref{#1}}

\newcommand{\lemref}[1]{Lemma~\ref{#1}}
\newcommand{\thmref}[1]{Theorem~\ref{#1}}
\newcommand{\propref}[1]{Proposition~\ref{#1}}

\newcommand{\itemref}[1]{\emph{\ref{#1}.}}

\newcommand{\PROOF}{\noindent{\it Proof : }} %%{\begin{proof}}
\newcommand{\ENDPROOF}{\hfill $\square$\newline}%%{\end{proof}}

\begin{document}

\title{Filled Julia sets of Chebyshev polynomials
%%\thanks{
%%The authors would like to thank
%%the Danish Council for Independent Research $|$ Natural Sciences for support via the
%%grant DFF -- 4181-00502.
%%The last author would also like to thank the Institute of Mathematical Sciences
%%of Stony Brook University for support and hosting during the %the {\color{red}{last}}
%%writings of the paper.}
}

\author[Petersen, Pedersen, Henriksen and Christiansen]
{Jacob Stordal Christiansen, Christian Henriksen, Henrik Laurberg Pedersen and Carsten Lunde Petersen}

\date{\today}

\maketitle

\begin{abstract}
We study the possible Hausdorff limits of the Julia sets and filled Julia sets of subsequences of 
the sequence of dual Chebyshev polynomials of a non-polar compact set $K\subset\C$ and compare such limits to $K$. 
Moreover, we prove that the measures of maximal entropy for the sequence of dual Chebyshev polynomials of $K$ 
converges weak* to the equilibrium measure on $K$.
\end{abstract}

\noindent {\em \small 2020 Mathematics Subject Classification: Primary: 42C05, Secondary: 37F10, 31A15}

\noindent {\em \small Keywords: Chebyshev polynomials, Julia set, Green's function}

\section{Introduction}%% and main results}
\label{intro}
%In this paper, we study the filled Julia sets of Chebyshev polynomials. 
Just like orthogonal polynomials, Chebyshev polynomials are central objects in numerical analysis, and yet very little seems to be known about their dynamical properties.
In this paper, we study the limit behavior of the filled Julia sets of dual Chebyshev polynomials and obtain results similar to those known for orthogonal polynomials (see \cite{CHPP} and \cite{Petersen-Uhre}).

Let $K\subset\C$ be an infinite compact set. 
The Chebyshev polynomial for $K$ of degree $n$ is the unique monic polynomial 
$T_n$ of minimal supremum-norm $||T_n||_{K,\infty}$. 
We denote by $\Om$ the unbounded connected component of $\C\sm K$
and define $J$ as the outer boundary of $K$, that is, $J=\partial\Om \subset \partial K\subset K$. 
Moreover, we denote by $\Po(K) := \C\Sm\Om$ the filled-in $K$, also known as the polynomial convex hull of $K$, and by $\Co(K)$ 
the Euclidean convex hull of $K$. 

We shall work with the \emph{dual Chebyshev polynomials}, that is, the polynomials $\T_n$ given by
$$\T_n(z) = T_n(z)/||T_n||_{K,\infty}$$ 
and write 
$$
\T_n(z) = \ga_n z^n + \textrm{lower order terms},
$$
where $\ga_n = 1/||T_n||_{K,\infty}$. Note that $\T_n$ is the unique degree $n$ polynomial
which is bounded by $1$ on $K$ and for which the leading coefficient is positive and maximal amongst all such polynomials.
In other words, these polynomials solve the maximization problem dual to the minimization problem of $T_n$ (see \cite{CSZ5}).

Denote by $\KK$ the set of non-empty compact subsets of $\C$. 
For $\{K_n\}_n\subset\KK$ a uniformly bounded sequence of compact sets, 
%%that is a sequence for which there exists $R>0$
%%such that $K_n\subset\D(R)$ for all $n$, 
we define compact limit sets 
 $\displaystyle{I :=\liminf_{n\to\infty} K_n}$ and $\displaystyle{S := \limsup_{n\to\infty} K_n}$, with $I\subset S$, by
\begin{align}
\liminf_{n\to\infty} K_n &:= \{ z\in\C \,|\,
\exists \, \{z_n\},\, K_n\ni z_n \underset{n\to\infty}\longrightarrow z\},\label{Hliminf}\\ 
\limsup_{n\to\infty} K_n &:= \{ z\in\C \,|\, \exists \, \{n_k\}, \, n_k\nearrow \infty
\textrm{ and }
\exists \, \{z_{n_k}\}, \, K_{n_k}\ni z_{n_k} \underset{k\to\infty}\longrightarrow z\}\label{Hlimsup}.
\end{align} 
The set $I$ may be empty, whereas $S$ is always non-empty. 
We will say that the sequence $\{K_n\}_n$ converges to $K$ and write $\lim_{n\to\infty} K_n = K$ if and only if 
$$I=S=K.$$ 
This defines a metrizable topology called the Hausdorff topology on $\KK$. 
For a brief introduction to the underlying complete metric on $\KK$, see Section~\ref{Hausdorffdistance}.

In \cite{Douady}, Douady proved that if $\{P_n\}_n$ is a sequence of polynomials of fixed degree $k\geq 2$ and if $P_n\to P$, then
$$
J \subset \liminf_{n\to\infty} J_n
\quad\textrm{and}\quad
\limsup_{n\to\infty} K_n \subset K,
$$
where $J_n$, $J$ are the Julia sets and $K_n$, $K$ the filled Julia sets for $P_n$ and $P$. 
In \cite{CHPP}, we proved potential theoretical analogs of these relations for the sequence of orthonormal polynomials 
defined by a Borel probability measure with non-polar compact support in $\C$.
These results were subsequently complemented with convergence statements for the sequence of measures of maximal entropy 
by Petersen and Uhre \cite{Petersen-Uhre}. 

In this paper, we prove versions of these theorems in the setting of Chebyshev polynomials. 
%As for notation, we denote by $\Om_n$ the attracted basin of $\infty$ for $\T_n$,
%by $K_n = \C\sm \Om_n$ the filled Julia set, and by $J_n =\partial K_n = \partial\Om_n$ the Julia set. 
%For $n\geq 2$, the common equilibrium measure $\omega_n$ for $K_n$ and $J_n$ 
%is also the unique measure of maximal entropy for $\T_n$ (see, e.g., \cite{Brolin}). 
%Note that in the case at hand, %the Chebyshev polynomials for a non-polar compact set 
%the inclusions are reversed.

\REFTHM{Main}
Let $K\subset\C$ be a non-polar compact set and let $\{\T_n\}_n$ 
be the associated sequence of dual Chebyshev polynomials.  
Then the corresponding sequence of filled Julia sets $\{K_n\}_{n\geq 2}$ 
is pre-compact in $\KK$ and for any limit point $K_\infty$ 
of a convergent subsequence $\{K_{n_k}\}_k$, we have that
$$
K\subset \Po(K_\infty) \subset
\Po(\limsup_{n\to\infty} K_n) \subset \Co(K).
$$
\ENDTHM
We conjecture --\;but are not able to prove -- that $K_\infty\sm K$ is a polar subset. 
From this it would follow that 
$\displaystyle{J = \partial \Po(K) \subset \liminf_{n\to\infty} J_n}$, which is similar to the orthogonal polynomial case.

\thmref{Main} gives upper and lower bounds on the possible limits of the filled Julia sets. The following theorem shows that in measure theoretical sense we do have convergence.
\REFTHM{Weakstarlmtsofomn}
Let $K\subset\C$ be a non-polar compact set, let $\T_n$ 
be the corresponding dual Chebyshev polynomials, and 
let $\om_n$ (for $n\geq 2$) denote the unique measure of maximal entropy for $\T_n$. 
Then the sequence $\{\om_n\}_n$ converges weak* to the equilibrium measure on $K$:
$$
\om_{n} \xrightarrow{\textrm{weak}^*} \om_K. 
$$
\ENDTHM

As for notation, we denote by $\Om_n$ the attracted basin of $\infty$ for $\T_n$,
by $K_n = \C\sm \Om_n$ the filled Julia set, and by $J_n =\partial K_n = \partial\Om_n$ the Julia set. 
Recall that for $n\geq 2$, the common equilibrium measure $\omega_n$ for $K_n$ and $J_n$ 
is also the unique measure of maximal entropy for $\T_n$ (see, e.g., \cite{Brolin}). 

\begin{figure}[h!] 
  \begin{center}
    \; \hfill \includegraphics[width=0.45\textwidth]{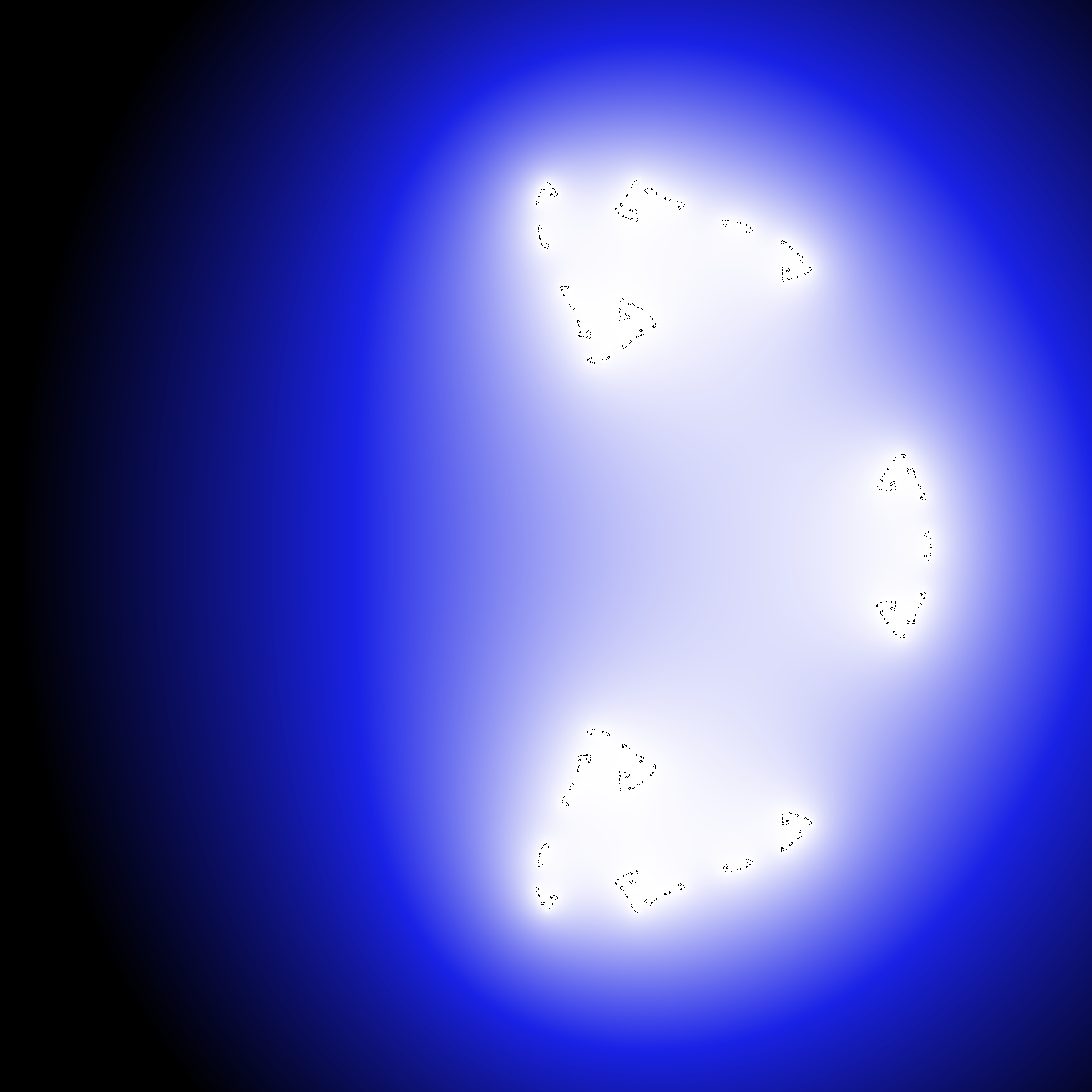}
    \hfill
    \includegraphics[width=0.45\textwidth]{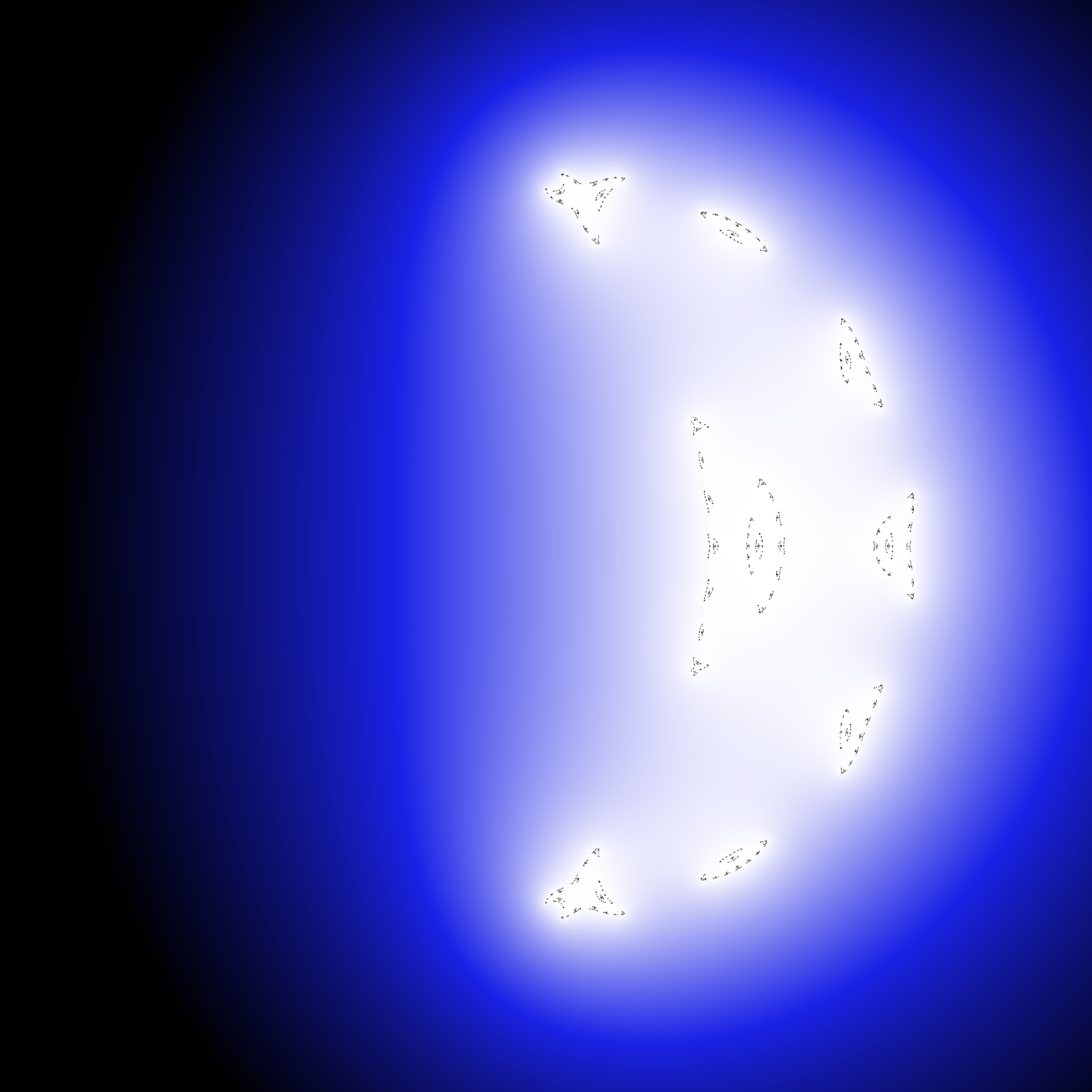} \hfill \;\\
    \vspace{0.02\textwidth}
    \;\hfill \includegraphics[width=0.45\textwidth]{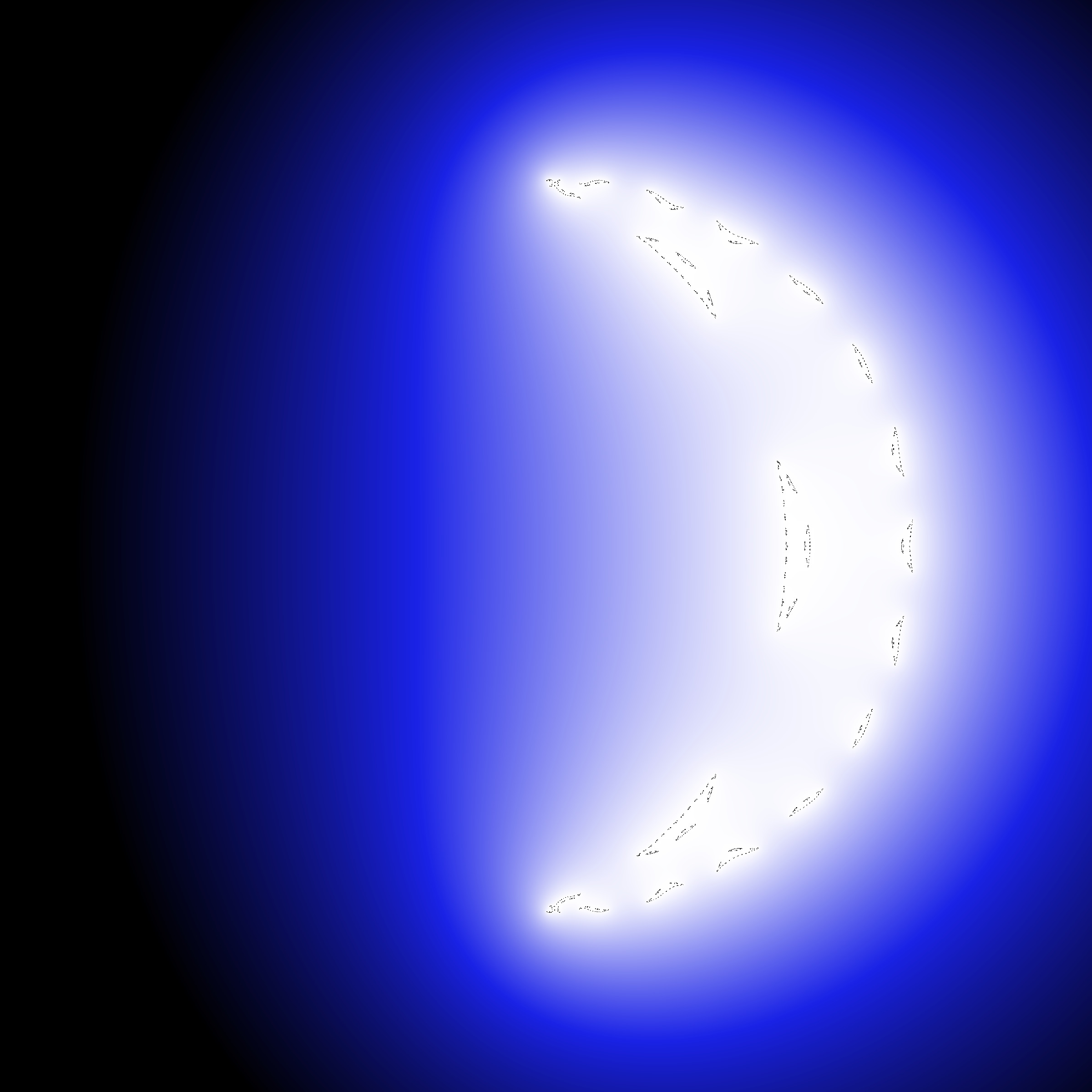}
    \hfill
    \includegraphics[width=0.45\textwidth]{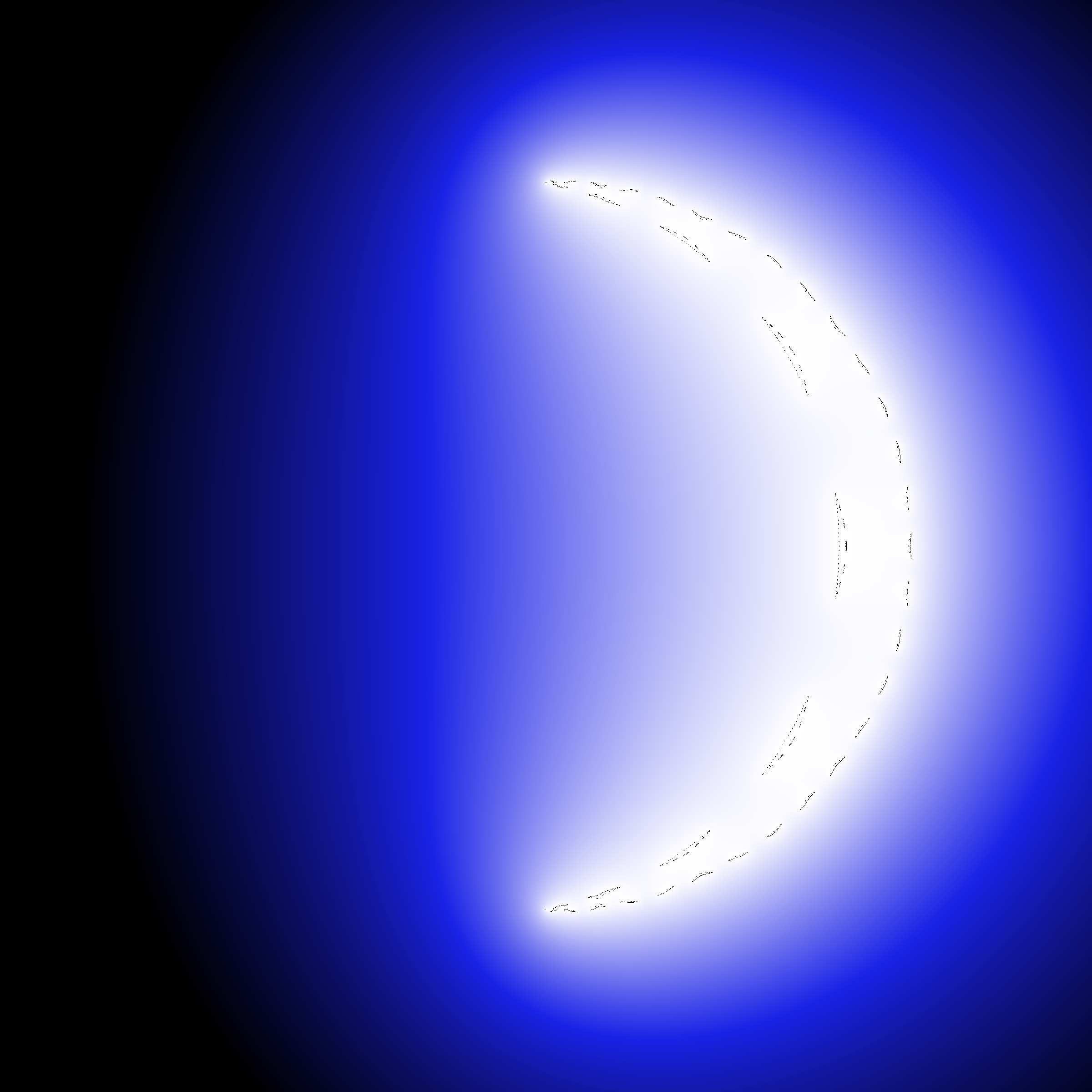} \hfill \; \\
  \end{center}
  \caption{The graphics illustrate the filled Julia sets $K_5$, $K_{10}$, $K_{20}$, and $K_{40}$ of numerically approximated Chebyshev polynomials associated to the arc $K = \{ e^{i\theta} : -\frac{\pi}{2} \leq \theta \leq \frac{\pi}{2} \}$.
  The windows are chosen so that the real and imaginary parts of $z$ ranges between $-1.5$ and $1.5$.
  The filled Julia set is drawn in black.
  The complement is drawn in blue, with a brightness decided by the Green's function associated to $K_n$, such that for values of $g$ close to zero, the color is light, whereas the color get increasingly darker as $g$ increases.
  The main theorem states that any accumulation point of a sequence $z_n \in K_n$ must be contained in the right half-disk, whereas the polynomial convex hull of any Hausdorff limit point $K_\infty$ of the sequence $\{K_n\}_n$ contains the arc $K$.}
    \label{fig:halfcircle}

\end{figure}

In figure \ref{fig:halfcircle}, we have pictured some filled Julia sets for the Chebyshev polynomials of $K = \{ e^{i\theta} : -\frac{\pi}{2} \leq \theta \leq \frac{\pi}{2} \}$ to illustrate the content of our main results.  See \cite{BE, Sch} for further details on the structure of these polynomials.

\section{Background}
We start by introducing the main players, namely potential theory, polynomial dynamics, Chebyshev polynomials, and Hausdorff distance.
\subsection{Potential theory in brief}
We shall only be concerned with the case where $K$ is non-polar. 
The outer boundary $J$ of a non-polar compact set $K$ supports 
a unique harmonic measure $\om_K$ for $\Om$ with respect to $\infty$ (on the Riemann sphere $\Chat$). 
This measure equals the so-called equilibrium measure of both $J$ and $K$, 
that is, the unique Borel probability measure supported on $K$ of maximal energy (following the sign convention of \cite{Ransford}).  
Indeed, defining the energy of a compactly supported Borel probability measure $\nu$ on $\C$ 
as the extended real number
\[
I(\nu) = \int\!\!\!\int_{\C\times\C} \log|z-w|\dnu(z)\dnu(w) \in [-\infty,\infty),
\]
then the maximal energy, among all Borel probability measures supported on $K$, exists and is realized by $\omega_K$.
The maximal energy $I(K) := I(\om_K)$ is finite by the non-polar hypothesis on $K$.

We denote by $g_\Om$ the Green's function for $\Om$ with pole at $\infty$, that is, the function 
\begin{equation}\label{eq:green}
g_\Om(z) = \int_\C\log|z-w|\d\om_K(w) - I(\om_K).
\end{equation}
The equillibrium measure $\om_K$ is also the distributional Laplacian $\Laplace g_\Om$ of the Green's function.
The capacity of $K$ is defined by $\Cpct(K) = \exp(I(K))$ and $g_\Om$ is the unique non-negative subharmonic function 
which is harmonic and positive on $\Om$, satisfies
\begin{equation}\label{eq:greenimplicit}
g_\Om(z) = \log|z| -\log\Cpct(K) + o(1) \;\mbox{ at infinity},
\end{equation}
and which is zero precisely on $K\sm E$,
where $E$ is the exceptional set defined by
\begin{equation*}
E = \{ z \in J \,|\, z \textrm{ is \emph{not} a Dirichlet regular boundary point}\} \subset J.
\end{equation*}
The set $E$ is an $F_\si$ polar set, see \cite[Theorems 4.2.5 and 4.4.9]{Ransford}.

\subsection{Polynomial dynamics in brief} \label{subsec:2.2}

For any polynomial $P$ of degree $d>1$, there exists $R>0$ such that $|P(z)| >|z|$ for $|z|\geq R$. 
It follows that for $|z|>R$ the orbit of $z$ under iteration by $P$ converges to $\infty$.
The basin of attraction of $\infty$ for $P$, denoted $\Om_P$, can therefore be written as
\begin{equation}
\label{basinofinfty}
 \Om_P = \{z\in\C \,|\, P^k(z) \underset{k\to\infty}\longrightarrow \infty \}
= \bigcup_{k\geq 0} P^{-k}(\C\sm\overline{\D(R)}).
\end{equation}
Here $P^k = \overset{k \textrm{ times}}{\overbrace{P\circ P \circ \ldots \circ P}}$, whereas $P^{-k}$ denotes the inverse image and $\D(R)$ is the open ball of radius $R$ centered at $0$.
It follows immediately that $\Om_P$ is open and completely invariant, that is,
$$P^{-1}(\Om_P) = \Om_P = P(\Om_P).$$
Denote by $K_P = \C\sm\Om_P\subseteq\overline{\D(R)}$ the filled Julia set for $P$ and
by $J_P = \partial\Om_P = \partial K_P$ the Julia set for $P$.
Then $K_P$ and $J_P$ are compact and also completely invariant.
Clearly, any periodic point (i.e., a solution of the equation $P^k(z) = z$ for some $k\in\N$) belongs to $K_P$,
so that $K_P$ is non-empty.
It follows from \eqref{basinofinfty} that the filled Julia set $K_P$ can also be described as the nested intersection
\begin{equation}
\label{filledJuliaset}
K_P = \bigcap_{k\geq 0} P^{-k}(\overline{\D(R)}).
\end{equation}
To ease notation, we denote the Green's function for $\Om_P$ with pole at infinity by $g_P$  (in lieu of the more cumbersome $g_{\Om_P}$).
It follows from \eqref{filledJuliaset} that $g_P$ satisfies
\EQN{Green_of_poly}
g_P(z) = \lim_{k\to\infty} \frac{1}{d^k}\log^+(|P^k(z)|/R) = \lim_{k\to\infty} \frac{1}{d^k}\log^+|P^k(z)|.
\ENDEQN
Here $\log^+$ is the positive part of $\log$ (i.e.,
$\log^+ x = \max \{\log x, 0\}$ for $x \geq 0$).
Thus $g_P$ vanishes precisely on $K_P$ and hence (\cite[Theorem 4.4.9]{Ransford}) every point in $J_P$ is a Dirichlet regular boundary point of $\Om_P$.
Moreover, denoting the leading coefficient of $P$ by $\gamma$, we see that
\EQN{g_P-recursion}
g_P(P(z)) = d\cdot g_P(z)
\quad\textrm{and}\quad
\Cpct(K_P) = \frac{1}{|\ga|^{\frac{1}{d-1}}}.
\ENDEQN
\subsection{Chebyshev polynomials in brief} \label{subsec:2.3}

%%Let $K\subset\C$ be any compact non-polar subset with $\Om$ the unbounded 
%%component of the complement and $g_\Om$ the corresponding Green's function. 
Recall that $K_n$ is the filled Julia set of the dual Chebyshev polynomial $\T_n$.
We have that
%\EQN{capicityconvergence}
\[
  \lim_{n\to\infty}\Cpct(K_n) = \lim_{n\to\infty}\frac{1}{\ga_n^{\frac{1}{n-1}}} 
  = \lim_{n\to\infty}\frac{1}{\ga_n^{\frac{1}{n}}} 
  %= \lim_{n\to\infty}||T_n||_{K,\infty}^{\frac{1}{n}}
   = \Cpct(K),
\]
%\ENDEQN
where the last equality sign is a classical result of Szeg\H{o} \cite{FaberFeketeSzegeo}, see also \eqref{nthrootreg} below.
Equivalently, 
\EQN{energyconvergence}
\lim_{n\to\infty}I(K_n) = I(K).
\ENDEQN

The Bernstein--Walsh Lemma (see, e.g., \cite[Theorem 5.5.7]{Ransford}) states that
\[
\forall\,n\in\N, \, \forall\,z\in\C : \frac{1}{n}\log |\T_n(z)| \leq g_\Om(z).
\]
Note that since $g_\Om \geq 0$, it follows trivially that also
\EQN{B-W-inq}
\forall\,n\in\N, \, \forall\,z\in\C : \frac{1}{n}\log^+ |\T_n(z)| \leq g_\Om(z).
\ENDEQN
The sequence $\T_n$  is $n$-th root regular, that is
\EQN{nthrootreg}
\lim_{n\to\infty} \frac{1}{n}\log|\T_n(z)| = g_\Om(z)
\ENDEQN
locally uniformly on $\C\sm\Co(K)$, see \cite[Theorem 3.9 (Chapter III)]{SaffTotik} or \cite[Theorem 3.2]{CSZ1}.
% Here locally uniformly is in the sense 
% of potential theory, that is, for a function $f$ and a sequence of functions $\{f_n\}_n$ on a subset $U\subset\C$ the statement
%$\lim_{n\to\infty} f_n(z) = f(z)$ locally uniformly on $U$ means
%\[
%\forall\,z\in U, \, \forall\,\{z_n\}_n\subset U :
 %\Bigl( z_n\underset{n\to\infty}\longrightarrow z \Longrightarrow  \lim_{n\to\infty} f_n(z_n) = f(z)\Bigr). 
%\]
%Note that for continuous functions on a compact subset $L\subset U$, 
%this is equivalent to locally uniform convergence in the classical sense.
%A similar terminology will be used for $\limsup_{n\to\infty} f_n$ and $\liminf_{n\to\infty} f_n$.

Finally, recall that for every $n$ the zeros of $\T_n$ are located in $\Co(K)$ (see \cite{Fejer} or, e.g., \cite[proof of Lemma 4]{Widom}).

%As a final preparation, we now briefly discuss the Hausdorff distance and Cara\-theo\-dory kernel convergence. 
\subsection{Hausdorff distance in brief}
\label{Hausdorffdistance}
Denote by $\KK$ the set of non-empty compact subsets of $\C$. 
The Hausdorff distance on $\KK$ is the natural choice in dynamical systems (see, e.g., \cite{Douady}) and we shall also use it here. 
Let us briefly recall the main definitions.
For $L, M \in \KK$, the Hausdorff semi-distance from $L$ to $M$ is given by
\begin{equation*}
\dH(L, M) := \sup\{\d(z, M) \,|\, z\in L\}
= \sup_{z\in L}\;\inf_{w\in M}\; |z-w|
\end{equation*}
and the Hausdorff distance between the two sets is given by 
\begin{equation*}
\DH(L, M) := \max\{\dH(L, M), \dH(M, L)\}.
\end{equation*}
The pair $(\KK, \dH)$ is a complete metric space.
A bounded sequence $\{K_n\}_n\subset\KK$ of compact sets is convergent if and only if 
$$\liminf_{n\to\infty} K_n = \limsup_{n\to\infty} K_n,$$
where $\liminf$ and $\limsup$ are defined as in \eqref{Hliminf} and \eqref{Hlimsup}. 
Moreover, any bounded sequence $\{K_n\}_n\subset \KK$ is sequentially pre-compact, 
that is, any subsequence has a convergent sub-subsequence. 
For compact subsets of the Riemann sphere, $\Chat$, we use instead the spherical metric. 
For compact subsets of $\C$ the spherical and the Euclidean metrics induce the same topology on $\KK$.

\section{Filled Julia sets and Green's functions for guided sequences of polynomials}
\label{section:3}
%The proofs in this section are similar to the proofs of similar statements about orthogonal polynomials. 
We start by introducing the notion of a \emph{guided} sequence of polynomials.
\DEF
Let $K$ be a non-polar compact set.
%with $\Om\subset\C$ the unbounded component of the complement and 
%$g_\Om$ the Green's function for $\Om$. 
We define a polynomial sequence of the form 
$\{P_n = \delta_nz^n+\textrm{lower order terms}\}_n$ to be $K$-guided if 
\begin{enumerate}
  \item the set $Z$ consisting of all zeros of all $P_n$ is bounded, and
  \item there exist constants $a>0$ and $b \in \R$ such that
\end{enumerate}
\EQN{lowerlogbound}
\liminf_{n\to\infty} \frac{1}{n}\log^+|P_n(z)| \geq a\cdot g_\Om(z) - b
\ENDEQN
locally uniformly on $\C\sm\Co(K)$, where $g_\Om$ is the Green's function for the unbounded component of the complement of $K$. 
\ENDDEF

\REFEXA{Chebychevdoit}
For a non-polar compact set $K$, the sequence of dual Chebyshev polynomials $\{\T_n\}_n$ is a $K$-guided sequence with $a=1$ and $b=0$.
This follows from $n$-th root regularity \eqref{nthrootreg} %of the sequence of Chebyshev polynomials 
and the fact that $Z\subset\Co(K)$. %according to \cite[proof of Lemma 4]{Widom}.
\ENDEXA

We now proceed with a general result which implies that there is a uniform upper bound on the size of the filled Julia sets $K_n$.

\REFPROP{basic_bound}
Let $K$ be a non-polar compact set and let $\{P_n\}_n$ be a $K$-guided sequence of polynomials.
Then there exists $R>0$ and $N\in\N$ such that
\EQN{containment}
\forall\, n\geq N : \quad K_n \subset P_n^{-1}(\overline{\D(R)}) \subset \D(R).%,
\ENDEQN
\ENDPROP
This statement is a vast generalisation of \cite[Lemma 2.1]{CHPP}.  Although the proof is similar, 
we provide it here for completeness.

\medskip

\PROOF
Fix $R>0$ such that $Z\subset \D(R)$ and $2\eps = \min\{ag_\Om(z)-b\;|\; |z| = R\}>0$. 
Then $C(0,R)$ is a compact subset of $\C\sm K$ and, in fact, $Z\cup K \subset \D(R)$. 
Hence, by the hypothesis \eqref{lowerlogbound}, we can choose $N_1$ such that 
\EQN{lowernormbound}
\forall\,n\geq N_1,\, \forall\,z \in C(0,R) : \frac{1}{n}\log|P_n(z)| \geq \eps.
\ENDEQN
Since all zeros of each $P_n$ is contained in $Z\subset\D(R)$, 
it follows from the minimum principle for harmonic functions that 
\eqref{lowernormbound} holds for $|z|\geq R$.
Choose any $N\geq N_1$ such that $R< \e^{N\eps}$. 
Then
$$
\forall\,n\geq N, \, \forall\,z, \, |z|\geq R : |P_n(z)| \geq  \e^{N\eps} > R.
$$
Hence $P_n^{-1}(\overline{\D(R)})\subset \D(R)$ for all $n\geq N$.
The inclusion $K_n \subset P_n^{-1}(\overline{\D(R)})$
follows immediately from \eqref{filledJuliaset}. The proof is complete.
\ENDPROOF

\REFCOR{limsupKnupperBound}
Let $K\subset\C$ be a non-polar compact set and $\{\T_n\}_n$ the sequence of dual Chebyshev polynomials for $K$. 
Then the corresponding sequence of filled Julia sets $\{K_n\}_{n\geq 2}$ is bounded 
and hence sequentially pre-compact in $\KK$. Moreover, 
\[
\limsup_{n\to\infty} K_n \subset \Co(K).
\]
\ENDCOR
\PROOF
According to \exaref{Chebychevdoit} and \propref{basic_bound}, there exists $R>0$ and $N$ such that 
\eqref{containment} is satisfied. 
Thus the sequence $\{K_n\}_n$ is bounded in $\C$ and hence pre-compact in $\KK$. 

Let $\Delta$ be any bounded topological disk with $\overline{\Delta}\subset\Chat\Sm\Co(K)$. 
Then the boundary $\partial\Delta$ is a compact subset of $\C\Sm\Co(K)$. 
Hence $2\eps := \min\{g_\Om(z)\;|\;z\in\partial \Delta\}>0$. 
%According to \eqref{nthrootreg} the hypothesis \eqref{lowerlogbound} is satisfied with $a=1$ and $b=0$ and 
%moreover the set $Z$ consisting of all zeros of all $T_n$ is contained in $\Co(K)$ and so is bounded. Hence the hypotheses of \propref{basic_bound} are satisfied.
%Let $R>0$ and $N$ be as in the conclusion of \propref{basic_bound}. 
Arguing as in the proof of \propref{basic_bound}, 
we now find that 
\[%\EQN{lowerlognormbound}
\forall\,n\geq N, \, \forall\,z \in \partial\Delta : \frac{1}{n}\log|\T_n(z)| \geq \eps.
\]%\ENDEQN
Consequently,
\EQN{escape}
\forall\,n\geq N, \, \forall\,z\in\partial\Delta : |\T_n(z)| \geq  \e^{N\eps} > R.
\ENDEQN
Since $Z\subset\Co(K)$, we have $Z\cap\overline{\Delta}=\emptyset$ and 
therefore \eqref{escape} holds on $\overline{\Delta}$ by the minimum principle for 
non-vanishing holomorphic functions. 
Moreover, $\overline\Delta\cap K_n=\emptyset$ for $n\geq N$ 
by \eqref{containment}.
Since $\Delta$ was arbitrary, we have $\limsup_{n\to\infty} K_n \subset \Co(K)$. 
\ENDPROOF
%This gives the upper inclusion of \thmref{Main}.

\REFPROP{g_comparison}
Let $0<C\leq 1\leq R$.
Then for any polynomial $P$ of degree $n\geq 2$ with filled Julia set $K$ satisfying  
$\Cpct(K)\geq C$ and 
$$
K\subset P^{-1}(\overline{\D(R)})\subset \D(R),
$$
the Green's function $g=g_P$ satisfies
\[%\EQN{logPnvsgnformula}
\Bigl\Vert \,g(z) - \frac{1}{n}\log^+|P(z)|\, \Bigr\Vert_\infty \leq \frac{M}{n},
\]%\ENDEQN
where $M = \log\frac{4R}{C}$.
\ENDPROP
This generalizes \cite[Prop.\;2.3]{CHPP} in that there is no reference to any particular sequence of polynomials 
nor extremality property. 
%We provide the details for completeness although it is similar.
%Note however that the statement here is completely free from the extremal polynomial setting. 
%There is even no mention of a sequence of polynomials og base compact set.

\medskip

\PROOF
Evidently, $\Cpct(K)<R$ by the general properties of capacity so that 
\[
\lvert \log\Cpct(K)\rvert \leq  \log\frac{R}{C}.
\]
By \eqref{eq:green} and \eqref{eq:greenimplicit}, the Green's functions $g$ 
can be written as
\[%\EQN{GreenInt}
g(z) = \log|z| - \log\Cpct(K) + \int_K \log|1-w/z| \,d\om(w),
\]%\ENDEQN
where $\om$ is the equilibrium measure on $K$.  

For $|z|\geq 2R$ and $w\in K$, we have $|w/z| < 1/2$ so that $\forall\,z, \, |z| \geq 2R$ : 
\EQN{g_log_comp}
\left| g(z) - \log\lvert z\rvert\right| \leq \lvert \log\Cpct(K)\rvert + \log 2 \leq \log\frac{R}{C} +\log 2 < M.
\ENDEQN
We divide the complex plane into the set $A = \{z \,|\, |P(z)| \leq 2R\}$ and its complement and 
estimate $|g(z) - \frac{1}{n}\log^+|P(z)||$ separately on each set.

For all $z\in\C\sm A$, we have $|P(z)| > 2R\geq 2$
so that $\log^+|P(z)| = \log|P(z)|$ and
$$
\left|g(z) - \frac{1}{n}\log^+|P(z)| \right| =
\left|\frac{1}{n}g(P(z)) - \frac{1}{n}\log|P(z)| \right| <\frac{M}{n}.
$$
Here we have used \eqref{g_P-recursion} and \eqref{g_log_comp} applied to $P(z)$.
For any $z\in \partial A$, we have $|P(z)|=2R$ so that
\begin{align*}
0 \leq g(z) &= \frac{1}{n}g(P(z))\\
 &= 
\frac{1}{n}\left( \log(2R) - \log\Cpct(K) +  \int_K \log|1-w/P(z)| \,d\om(w)\right)\\
&< \frac{1}{n}\left(\log\frac{2R}{C} + \log 2\right) = \frac{M}{n}.
\end{align*}
Hence, by the maximum principle for subharmonic functions,
$g(z) < {M}/{n}$ on $A$.
Similarly, 
$$
0 \leq \frac{1}{n}\log^+|P(z)| \leq \frac{\log 2R}{n} < \frac{M}{n}
$$
on $A$ by construction. Finally,
$$
\left| g(z) - \frac{1}{n}\log^+|P(z)| \right| \leq \max\{g(z), \frac{1}{n}\log^+|P(z)|\} \leq \frac{M}{n}
$$
on $A$. This completes the proof.
\ENDPROOF

The above proposition gives us uniform bounds on the Green's functions $g_n$ for the dual Chebyshev polynomials 
$\T_n$ (cf. \eqref{Green_of_poly}). 
\REFCOR{Chebychev_g_comparison}
For a non-polar compact set $K$ and the corresponding dual Chebyshev polynomials $\T_n$ 
with Green's functions $g_n$, there exists $N\in\N$ and $M>0$ such that
\begin{equation}\label{logPnvsgnformulaB}
\forall\, n\geq N :\quad \Bigl\Vert \,g_n(z) - \frac{1}{n}\log^+|\T_n(z)|\, \Bigr\Vert_\infty \leq \frac{M}{n}.
\end{equation}
\ENDCOR
\PROOF
By \exaref{Chebychevdoit}, there exists $R>1$ and $N\in\N$ such that the filled Julia sets $K_n$ for $\T_n$ satisfy \eqref{containment}. 
Moreover, $\Cpct(K_n)\to\Cpct(K) > 0$ as $n\to\infty$ 
so there exists $C\in(0, 1)$ such that $C < \Cpct(K_n)$ for all $n$. 
Thus \eqref{logPnvsgnformulaB} follows by applying \propref{g_comparison} to each $\T_n$. 
%with $M = \log\frac{4R}{C}$.
\ENDPROOF

\REFCOR{boundinggnandK}
For a non-polar compact set $K$ and the corresponding dual Chebyshev polynomials $\T_n$ with Green's functions $g_n$,
there exists $N\in\N$ and $M>0$ such that
\EQN{limsup}
\forall\,n\geq N, \, \forall\; z\in\C : \quad g_n(z) \leq g_\Om(z) + {M}/{n}.
\ENDEQN
In particular,
\[
\limsup_{n\to\infty} g_n(z) \leq g_\Om(z)
\]
uniformly on $\C$. 
Moreover, 
\[
\lim_{n\to\infty} g_n(z) = g_\Om(z)
\]
locally uniformly on $\C\sm\Co(K)$ and for every $n\geq N$, we have
\EQN{Moverncontainment}
K \subseteq \{z\;|\; g_n(z) \leq {M}/{n}\}.
\ENDEQN
\ENDCOR
\PROOF
As for \eqref{limsup}, combine the Bernstein--Walsh Lemma in the form \eqref{B-W-inq} with the above Corollary. 
Then the $\limsup$ statement immediately follows. 
Similarly, the $\liminf$ statement follows from the $n$-th root regularity of the Chebyshev polynomials and the above Corollary. 
As for \eqref{Moverncontainment}, combine the Corollary with 
$||\T_n||_{K,\infty} =1$ so that $\frac{1}{n}\log^+|\T_n(z)| \equiv 0$ on $K$.
\ENDPROOF

\section{Proofs of the main theorems}
We are now ready to present the proofs of our main results.
\subsection{Proof of \thmref{Main}.}
%\subsection{Caratheodory kernel convergence}
%The notion of Caratheodory kernel convergence makes sense in any metric space. 
%We shall use it here in the Riemann sphere  equipped with the spherical metric. 
The notion of Caratheodory convergence is the last ingredient we need to prove \thmref{Main}.
Recall that a pointed domain is a pair $(U, z)$, where $z\in U\subset \Chat$ and $U$ is a domain 
(i.e., connected and open).
\REFDEF{def:caratheodory}
Let $(U_n, z_n)_n$ be a sequence of pointed domains. 
We say that the sequence converges to the pointed domain $(U_\infty, z_\infty)$ if
\ENUM
\item
$z_n \to z_\infty$, % as $n\to\infty$,
\item\label{CKctwo}
any compact subset $K \subset U_\infty$ is a subset of $U_n$ 
for all but finitely many $n$, and
\item
for any open connected set $V$ with $z_\infty\in V$, 
if $V\subset U_n$ for infinitely many $n$, then $V\subset U_\infty$.
\ENDENUM
\ENDDEF

\REM
There is an alternative characterization of Caratheodory convergence.
The pointed domains $(U_n, z_n)_n$ converge to $(U_\infty, z_\infty)$ if and only if the following two properties hold:
\ENUMi
\item$z_n\to z_\infty$, and 
\item
for any compact subset $K\subset\Chat\Sm\{z_\infty\}$ and for any subsequence $(n_k)_k$, 
if $\Chat\Sm U_{n_k} \longrightarrow K\subset\Chat\Sm\{z_\infty\}$ 
in the Hausdorff distance on compact subsets of $\Chat$, %for some compact subset $K\subset\Chat\Sm\{z_\infty\}$, 
then the domain $U_\infty$ is equal to the component of $\Chat \Sm K$ containing $z_\infty$.
\ENDENUM
\ENDREM

Using the formulation of Caratheodory convergence as in the remark, one can show that even if $z_n\to z_\infty$, a pointed sequence $(U_n, z_n)_n$ may fail to converge either because (for some subsequence) $\d(z_n,\partial U_n) \to 0$ or because there are at least two subsequences $U_{n_k}$ and $U_{n_m}$ such that $\Chat\Sm U_{n_k}$ and $\Chat\Sm U_{n_m}$ converge to $K_1$ and $K_2$, respectively, but the connected components of $\Chat\Sm K_1$ and $\Chat\Sm K_2$ containing $z_\infty$ are different. 
However, if $\d(z_n,\partial U_n)$ is bounded uniformly from below by some $r>0$, 
then the sequence $(U_n, z_n)_n$ is sequentially pre-compact. 
In particular, any sequence of pointed domains $(\Om_n, \infty)$ in the Riemanns sphere 
with $K_n\subset \D(R)$ for some fixed $R>0$ is sequentially compact.

%%\LEM
%%If $\Om_n\subset\Chat$ is a sequence of domains with compact complements 
%%$K_n=\C\Sm\Om_n \subset\D(R)$ for some $R>0$. 
%%Then the sequence $(\Om_n,\infty)$ is sequentially pre-compact for the Caratheodory topology on pointed subsets of the Riemann spheres.
%%\ENDLEM
%%\PROOF
%%Given any subsequence $\{n_k\}_k$ of $\{n\}_n$, 
%%we can suppose possibly passing to a further subsequence also denoted 
%%$\{n_k\}_k$ that the compact complements $K_{n_k}$ converge to 
%%some compact set $K\subset\overline{\D(R)}$. Let $\Om$ denote the 
%%connected component of $\Chat\Sm K$ containing $\infty$. 
%%Then $(\Om_{n_k},\infty)$ converges to $(\Om,\infty)$ in the sense of Caratheodory.
%%\ENDPROOF
Given a sequence $(g_n)_n$ of Green's functions with pole at $\infty$ for domains $\Om_n\subset\Chat$ with compact complements $K_n=\C\Sm\Om_n$, there is no general relation between limits of subsequences $g_{n_k}$ and the question of Caratheodory convergence of the corresponding subsequence of pointed domains 
$(\Om_{n_k},\infty)$.
However, the following Proposition gives an upper bound 
on such Caratheodory limits.
\REFPROP{boundingK}
Let $(g_n)_n$ be a sequence of Green's functions with pole at $\infty$
for domains $\Om_n\subset\Chat$,
and let $K' = \{ z \mid g_n(z) \to 0 \}$.
Then 
$$
K'\cap\Om_\infty = \emptyset
$$
for any domain $\Om_\infty$ such that a subsequence $(\Om_{n_k},\infty)_k$ converges to $(\Om_\infty,\infty)$ 
in the sense of Caratheodory.
\ENDPROP
\PROOF
It suffices to consider the case where $(\Om_n,\infty) \longrightarrow (\Om_\infty,\infty)$ 
in the sense of Caratheodory. 
Let $z_0 \in \Om_\infty$ be arbitrary.
Choose a domain $U$, $\{z_0,\infty\}\subset U\subset \overline{U}\subset \Om_\infty$, having a Green's function $g_U$ with pole at $\infty$.
We may take $U$ to be a Jordan disk. 
In view of definition \ref{def:caratheodory} part \ref{CKctwo}, 
there exists $n_0$ such that $\overline{U} \subset \Om_n$ for all $n \geq n_0$.
By \cite[Cor.\;4.4.5]{Ransford}, $0 < g_U(z) \leq g_{\Omega_n}(z)$ for any $z\in U$ when $n\geq n_0$. 
In particular, $0< g_U(z_0) \leq g_n(z_0)$ for all $n\geq n_0$.
Hence $z \notin K'$.
\ENDPROOF

The upper bound in \thmref{Main} follows from \corref{limsupKnupperBound}. 
So it only remains to prove that
$$
\Om \supset \Om_\infty,
\quad \textrm{or equivalently,} \quad
K\cap \Om_\infty = \emptyset
$$
for any Caratheodory limit point $(\Om_\infty,\infty)$ 
of a convergent subsequence $(\Om_{n_k}, \infty)_k$.
However, according to \corref{boundinggnandK} the Green's functions 
$g_n$ for $\Om_n$ satisfies that $g_n(z) \to 0$ on $K$.
So the statement follows from \propref{boundingK}.

\subsection{Proof of \thmref{Weakstarlmtsofomn}.}
We need a few more auxillary results in order to prove this result. First of all, we shall use the following Lemma which is a refinement of \cite[Lemma 1.3.2]{StahlandTotik}. 
For a proof, the reader is referred to \cite[Lemma 3.3]{Petersen-Uhre}.
\REFLEM{deepbound}
Let $V, K\subset\C$ be compact sets with $V$ contained in the unbounded component of $\C \setminus K$, and let $b\in(0, 1)$ be arbitrary. 
Then there exists $M = M(b, V, K) \in\N$ such that for 
$M$ arbitrary points $x_1, x_2, \ldots ,x_M \in V$, there exists $M$ points 
$y_1, y_2, \ldots y_M\in \C$ for which the rational function 
\[
r(z) = \prod_{j=1}^M \frac{z-y_j}{z-x_j}
\]
has supremum norm on $K$ bounded by $b$, that is,
$
||r||_K \leq b.
$
\ENDLEM

We shall also make use of the following result.

\REFPROP{boundonnoofpreimages}
Let $K\subset\C$ be a non-polar compact set with corresponding dual Chebyshev polynomials $\T_n$.
For any $R>0$ and any compact set $V\subset \C$ with $V\cap \Po(K)=\emptyset$, there exists $M=M(K, R, V) \in\N$ 
such that for any $w$, $|w| \leq R$, and 
for any $n$, the number of pre-images of $w$ in $V$ 
under the polynomial $\T_n$ is less than $M$. 
In symbols,
$$
\# [ \T_n^{-1}(w) \cap V ] < M.
$$
\ENDPROP
\PROOF
Fix $R>0$ and any compact set $V\subset \C$ with $V\cap \Po(K)=\emptyset$. 
Let $b = 1/(1+R)$ and let $M=M(b, V, \Po(K))$ be as in \lemref{deepbound}. 
For each $n$, let $T_n(z)$ be the usual monic Chebyshev polynomial 
of degree $n$ for $K$. 
Then $T_n$ is the unique degree $n$ monic polynomial 
of minimal sup-norm $1/\ga_n$ on $K$. 

Suppose towards a contradiction that for some $w$, $|w|\leq R$, and some $n$, the 
equation $\T_n(z) = w$ has at least $M$ solutions $x_1, \ldots x_M \in V$. 
Let $r$ be the rational function given by \lemref{deepbound} such that 
$||r||_{K} \leq b$ and set $q(z) := r(z)\cdot(\T_n(z)-w)/\ga_n$. 
Then $q$ is a monic polynomial of degree $n$ and
\begin{align*}
  ||q||_K &\leq \frac{||r||_K}{\ga_n} ||\T_n(z)-w||_K \leq\frac{b}{\ga_n}\sqrt{1+|w|}\\
  &\leq\frac{\sqrt{1+R}}{1+R}\cdot ||T_n||_K < ||T_n||_K.
\end{align*}
This contradicts the fact that 
$T_n$ has minimal sup-norm on $K$ among all monic degree $n$ polynomials.
\ENDPROOF

\REFCOR{Knboundedness} 
Let $K\subset\C$ be a non-polar compact set and suppose $V\subset \C$ with $V\cap \Po(K)=\emptyset$. Then
there exists $M=M(K;V) \in\N$ such that for any of the dual Chebyshev polynomials $\T_n$ with $n\geq 2$ and  
any $z\in K_n$, the number of pre-images of $z$ in $V$ under $\T_n$ is less than $M$. 
%That is
%$$
%\#(T_n^{-1}(z) \cap V) < M.
%$$
\ENDCOR
\PROOF
It follows from \propref{basic_bound} that there exists $R>0$ such that 
$$
K_n \subset \T_n^{-1}(\overline{\D(R)}) \subset \D(R)
$$ 
for all $n\geq 2$, where $K_n$ is the filled Julia set of $\T_n$. 
The Corollary follows immediately from this and \propref{boundonnoofpreimages}.
\ENDPROOF

\PROOF (of \thmref{Weakstarlmtsofomn})
According to \corref{limsupKnupperBound}, the sequence of equillibrium measures $\om_n$ for $K_n$ have uniformly bounded support and so the sequence of such measures is pre-compact for the weak-* topology.
Furthermore, according to Brolin \cite{Brolin} (see also Lyubich \cite{Lyubich}), 
the equilibriums measure $\om_n$ is also the unique invariant balanced measure for $\T_n$,  
that is, it is the unique probability measure $\om$ on $\C$ such that for any measurable function 
{\mapfromto f \C \C},
\[ %\EQN{balanced}
\int_{\C} f(z)d\om(z) = \frac{1}{n}\int_{\C} \left(\sum_{w, \T_n(w) = z} f(w)\right) d\om(z).
\] % \ENDEQN
Let $V\subset\C$ be a compact subset with $V\cap \Po(K) = \emptyset$ 
and let $M\in\N$ be as in \corref{Knboundedness}. 
Then for the measureable function $1_V$ (i.e., the indicator function for $V$), we have 
\[
\om_n(V) = \int_ {\C}1_V(z)d\om_n(z) = \frac{1}{n}\int_{\C} \left(\sum_{w, \T_n(w) = z} 1_V(w)\right) d\om_n(z) 
\leq \frac{M}{n}\underset{n\to\infty}\longrightarrow 0.
\]
This proves that for any weak limit $\nu$ of a convergent subsequence $\{\om_{n_k}\}_k$, the support $S(\nu)$ is contained in $\Po(K)$. 
Furthermore, by \eqref{energyconvergence} and \cite[Lemma 3.3.3]{Ransford}), we have
\[
I(\nu) \geq \limsup_{k\to\infty} I(\om_{n_k}) = \lim_{n\to\infty} I(\om_n) =  I(\om_K) = I(\Po(K)).
\]
Hence $\nu=\om_K$ since $\om_K$ is the unique measure of maximal energy $I(K)$.
As the limit is unique, we in fact have 
\[\om_n \xrightarrow{\textrm{weak}^*} \om_K\]
and this proves  \thmref{Weakstarlmtsofomn}.
\ENDPROOF



%\authorrunning{Petersen, Pedersen, Henriksen and Christiansen} % if too long for running head

              Jacob Stordal Christiansen,
              Lund University, Centre for Mathematical Sciences, Box 118, 22100 Lund, Sweden\\
%              Tel.: +123-45-678910\\
%              Fax: +123-45-678910\\
              \texttt{stordal@maths.lth.se}           %  \\
\medskip

      Christian Henriksen,
              DTU Compute, Technical University of Denmark, Build. 303B, 2800 Kgs. Lyngby, Denmark\\
%              Tel.: +123-45-678910\\
%             Fax: +123-45-678910\\
              \texttt{chrh@dtu.dk}           %  \\

\medskip

      Henrik Laurberg Pedersen,
               Department of Mathematical Sciences, University of Copenhagen,
	       Universitetsparken 5, 2100 Copenhagen, Denmark\\
%              Tel.: +123-45-678910\\
%              Fax: +123-45-678910\\
              \texttt{henrikp@math.ku.dk}           %  \\

\medskip

      Carsten Lunde Petersen,
              Department of Science and Environment, Roskilde University, 4000 Roskilde, Denmark \\
%              Tel.: +45-46742076\\
%              Fax: +123-45-678910\\
              \texttt{lunde@ruc.dk}           %  \\

\end{document}